\documentclass[a4paper,11pt,leqno]{article}
\usepackage{amsmath, amssymb, amsthm}
\usepackage{graphicx} 

\usepackage[latin1]{inputenc}

\numberwithin{equation}{section}

\begin{document}

\def\FF{{\mathcal F}} \def\SS{{\mathcal S}} \def\C{{\mathbb C}}
\def\F{{\mathbb F}} \def\Z{{\mathbb Z}} \def\S{{\mathbb S}}
\def\LL{{\mathcal L}} \def\N{{\mathbb N}} \def\K{{\mathbb K}}
\def\SSS{{\mathfrak S}} \def\HSS{{\widehat{\mathcal S}}}
\def\bO{{\bf O}} \def\BB{{\mathcal B}}


\title{\bf{HOMFLY-PT skein module of singular links in the three-sphere}}

\author{\textsc{Luis Paris\footnote{Partially supported by the {\it Agence Nationale de la Recherche} ({\it projet Théorie de Garside}, ANR-08-BLAN-0269-03).} 
and Emmanuel Wagner\footnote{Partially supported by the {\it Agence Nationale de la Recherche} ({\it ANR Jeunes Chercheurs Vaskho}, ANR-11-JS01-002-01) and the {\it Conseil Régional de Bourgogne} ({\it Bourse Faber} X110CVHCP-2011)}
}}

\date{\today}

\maketitle

\begin{abstract}
\noindent
For a ring $R$, we denote by $R[\LL]$ the free $R$-module spanned by the isotopy classes of singular links in $\S^3$.
Given two invertible elements $x,t \in R$, the HOMFLY-PT skein module of singular links in $\S^3$ (relative to the triple $(R,t,x)$) is the quotient of $R[\LL]$ by local relations, called skein relations, that involve $t$ and $x$.
We compute the HOMFLY-PT skein module of singular links for any $R$ such that $(t^{-1}-t+x)$ and $(t^{-1}-t-x)$ are invertible.
In particular, we deduce the Conway skein module of singular links.
\end{abstract}

\noindent
{\bf AMS Subject Classification.} Primary: 57M25.


\section{Introduction}

HOMFLY-PT Skein modules, introduced by Przytycki and Turaev \cite{Przyt1,Turae1}, are $3$-dimensional quantum invariants which can be thought as a way to extend to three manifolds famous quantum link invariants such as the Jones polynomial \cite{Jones2,Jones1} and the HOMFLY-PT polynomial \cite{HOMFLY,Jones1,PrzTra1}. 
Given a three manifold $M$, a ring $R$, and two invertible elements $x,t \in R$, let $\FF(M)$ be the free $R$-module generated by the isotopy classes of (oriented) links in $M$. The HOMFLY-PT skein module $\SS(M)=\SS(M;R,t,x)$ is the quotient of $\FF(M)$ by local relations, called skein relations, that involve $t$ and $x$. They are computed in very few cases, see for examples \cite{Przyt2, GilZho1, GilZho2}.

\bigskip\noindent
Instead of trying to compute this invariant for different three-manifolds, one can also consider other sets of isotopy classes of topological objects in a given three-manifold. One can, for instance, consider embedded surfaces in a three-manifold \cite{AsaFro1}, or enlarge the class of links by considering graphs or singular links into a three-manifold.

\bigskip\noindent
In \cite{Paris1} the first author considered the HOMFLY-PT skein module of singular links in the three-sphere. Using heavily previous work by Rabenda and himself  \cite{ParRab1} on Markov traces on singular Hecke algebras, he was able to determine this skein module when $R = \C(t, x)$ is the two variable field of fractions over $\C$.

\bigskip\noindent
In this paper we determine the HOMFLY-PT skein module of singular links in the three-sphere for any ring $R$ such that $(t^{-1}-t+x)$ and $(t^{-1}-t-x)$ are invertible.
Our result concerns the ring $R=\C(t,x)$ already considered in \cite{Paris1}, but it also concerns other rings of coefficients such as $\F_q(t,x)$, where $\F_q$ is a finite field of cardinality $q$, or $\Z[x^{\pm 1}]$ (with $t=1$), corresponding to the so-called {\it Conway skein module}.
Moreover, the paper is self-contained, and the proofs are simpler and more direct than those given in \cite{ParRab1,Paris1} for the case $R=\C(x,t)$.

\bigskip\noindent
Before going into details of the statements and proofs, we wish to point out two distinguished points in the paper. 
As for the links in thickened surfaces (see \cite{Przyt2}), we can define the product of two links by placing one above the other, and this (commutative) multiplication equips the HOMFLY-PT skein module of singular links with an algebra structure.
However, unlike in the case of thickened surfaces \cite{Przyt2}, this multiplication can be trivial, thus the HOMFLY-PT skein module of singular links equipped with this operation is not necessarily a polynomial algebra.
In this paper we define another operation by means of connected sums, and we show that the HOMFLY-PT skein module of singular links equipped with this second operation is a polynomial algebra over $R$ in two variables.

\bigskip\noindent
The second distinguished point is that, in order to treat the case where $R$ is of characteristic different from zero, we introduce a new category of knot-like objects: the ``ordered singular links''.
These are simply singular links whose singular points are linearly ordered.
In order to determine the HOMFLY-PT skein module of singular links, we need first to determine the HOMFLY-PT skein module of ordered singular links.
This is a free algebra over $R$ freely generated by two elements.


\section{Definitions and statements}

Recall that a \emph{singular (oriented) link} with $n$ components in $\S^3$ is an immersion of $n$ circles in $\S^3$ which admits only finitely many singularities that are all transverse double points. As usual in knot theory, we reduce the study of objects in dimension three to the study of diagrams (generic projection) in the plane. By a result of Kauffman \cite{Kauff1}, one has a version of Reidemeister theorem for singular links. We will denote by $\LL$ the set of isotopy classes of singular links, and, for $d \in \N$, by $\LL_d$ the set of singular links with precisely $d$ singular points.

\bigskip\noindent
Let $R$ be a commutative ring, and let $t,x$ be two invertible elements in $R$. Let $R[\LL]$ denote the free $R$-module freely generated by $\LL$, and, for $d \in \N$, let $R[\LL_d]$ denote the free $R$-module freely generated by $\LL_d$. The \emph{HOMFLY-PT skein module} of singular links (relative to the triple $(R,t,x)$), denoted by $\SS(R,t,x)$, is defined to be the quotient of $R[\LL]$ by the relations
\begin{equation} \label{skein}
x\, L_0 = t^{-1}\, L_+ - t\,L_-\,,
\end{equation}
for all triples of singular links $(L_0,L_+,L_-)$ such that $L_0,L_+,L_-$ have the same diagram except in a disk where they are like in Figure 2.1. On the other hand, for $d \in \N$, we denote by $\SS_d (R,t,x)$ the quotient of $R[\LL_d]$ by the relations (\ref{skein}). Note that 
\[
\SS (R,t,x) = \bigoplus_{d=0}^\infty \SS_d(R,t,x)\,.
\]

\begin{figure}[tbh]
\bigskip\bigskip\centerline{
\setlength{\unitlength}{0.4cm}
\begin{picture}(10,1.2)
\put(0,0.4){\includegraphics[width=4cm]{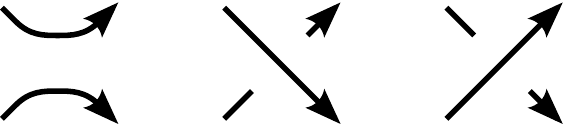}}
\put(0.5,0){\small $L_0$}
\put(4.5,0){\small $L_+$}
\put(8.5,0){\small $L_-$}
\end{picture}} 

\smallskip
\centerline{{\bf Figure 2.1.} The links $L_0$, $L_+$ and $L_-$.}
\end{figure}

\bigskip\noindent
As pointed out in the introduction, there are two ways to equip $\SS(R,t,x)$ with an algebra structure. 
For $L_1,L_2\in \LL$, we denote by $L_1 \cdot L_2$ the singular link obtained placing $L_1$ above $L_2$.
it is easily seen that this multiplication induces an operation on $\SS(R,t,x)$, still denoted by $\cdot$, and that $\SS(R,t,x)$ endowed with this operation is a graded algebra. 

\bigskip\noindent
{\bf Remark.}
Note that $(\SS(R,t,x), \cdot)$ is not a priori a unit algebra.
We could have added the emptyset to $\LL$ in order to have a unit, but, in that case, we would have needed to add the relation $(t^{-1}-t) \emptyset = x \bO$ in the definition of $\SS(R,t,x)$, where $\bO$ denotes the trivial knot. 
Moreover, if $t=\pm 1$, then the product of two elements is always zero (this will be explained in more detail later), so we do not like to have a unit in that case. 
Since, in this paper, we will mainly consider the operation defined below, and this other operation has a unit in $\SS(R,t,x)$, the trivial knot $\bO$, we do not assume the emptyset to lie in $\LL$.  

\bigskip\noindent
We turn now to define the second operation. A \emph{connected sum} of two singular links $L_1$ and $L_2$ is defined to be a singular link $L$ obtained from $L_1$ and $L_2$ like in Figure 2.2. 

\begin{figure}[tbh]
\bigskip\bigskip\centerline{
\setlength{\unitlength}{0.4cm}
\begin{picture}(13,3)
\put(0,0.5){\includegraphics[width=5.2cm]{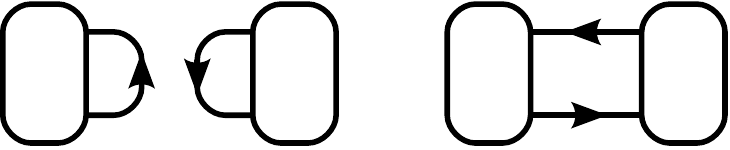}}
\put(0.4,1.6){\small $L_1$}
\put(4.8,1.6){\small $L_2$}
\put(8.3,1.6){\small $L_1$}
\put(11.7,1.6){\small $L_2$}
\put(10.2,0.2){\small $L$}
\end{picture}} 

\smallskip
\centerline{{\bf Figure 2.2.} Connected sum of $L_1$ and $L_2$.}
\end{figure}

\bigskip\noindent
If $L$ is a connected sum of two singular links $L_1$ and $L_2$, then $L$ and $L_1 \cdot L_2$ are related in $\SS (R,t,x)$ by the following equality.
\begin{equation}\label{deuxoperations}
(t^{-1}-t)L = x(L_1 \cdot L_2)\,.
\end{equation}
This equality is easily deduced from the graphical computation in $\SS(R,t,x)$ illustrated in Figure~2.3.

\begin{figure}[tbh]
\begin{gather*}
(t^{-1}-t)
\parbox[c]{2.3cm}
{\setlength{\unitlength}{0.3cm}
\begin{picture}(7,4)
\put(0,0){\includegraphics[width=2.1cm]{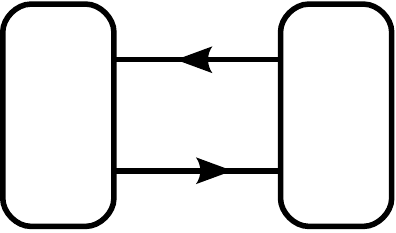}}
\put(0.5,1.6){\small $L_1$}
\put(5.5,1.6){\small $L_2$}
\end{picture}}
=t^{-1}
\parbox[c]{2.9cm}
{\setlength{\unitlength}{0.3cm}
\begin{picture}(9,4)
\put(0,0){\includegraphics[width=2.7cm]{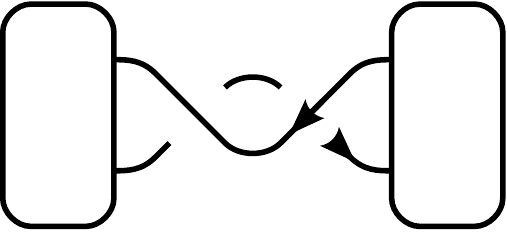}}
\put(0.5,1.6){\small $L_1$}
\put(7.5,1.6){\small $L_2$}
\end{picture}}
-t
\parbox[c]{2.9cm}
{\setlength{\unitlength}{0.3cm}
\begin{picture}(9,4)
\put(0,0){\includegraphics[width=2.7cm]{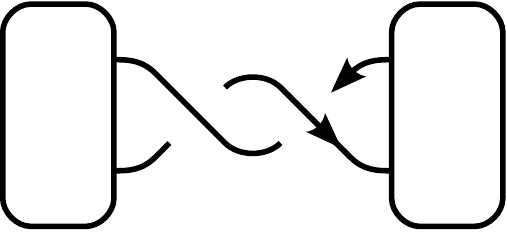}}
\put(0.5,1.6){\small $L_1$}
\put(7.5,1.6){\small $L_2$}
\end{picture}}\\
=x \left(
\parbox[c]{2.9cm}
{\setlength{\unitlength}{0.3cm}
\begin{picture}(9,4)
\put(0,0){\includegraphics[width=2.7cm]{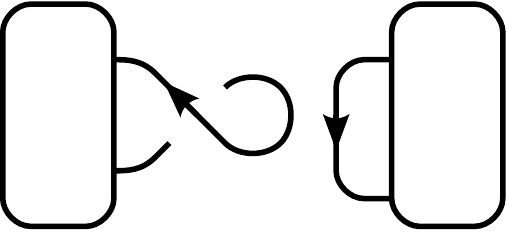}}
\put(0.5,1.6){\small $L_1$}
\put(7.5,1.6){\small $L_2$}
\end{picture}}\right)
=x \left(
\parbox[c]{2.3cm}
{\setlength{\unitlength}{0.3cm}
\begin{picture}(7,4)
\put(0,0){\includegraphics[width=2.1cm]{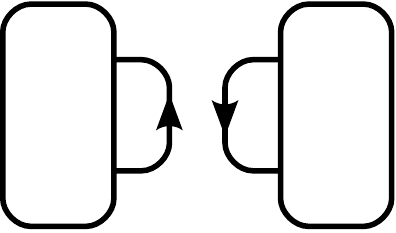}}
\put(0.5,1.6){\small $L_1$}
\put(5.5,1.6){\small $L_2$}
\end{picture}}\right)
\end{gather*}
\centerline{{\bf Figure 2.3.} Graphical computation.}
\end{figure}

\bigskip\noindent
A connected sum is not unique in general, but it is unique inside $\SS(R,t,x)$. This is easily seen if $R$ is an integral domain and $t \neq \pm 1$, because of Equality (2.2).
In the more general framework, this is still true but less obvious: 

\bigskip\noindent
{\bf Lemma 2.1.}
{\it Two connected sums of two singular links $L_1$ and $L_2$ represent the same element in $\SS(R,t,x)$.}

\bigskip\noindent
{\bf Proof.}
Lemma 2.1 is a straightforward consequence of the graphical computation in $\SS(R,t,x)$ shown in Figure 2.4.
\qed

\begin{figure}[tbh]
\begin{gather*}
\parbox[c]{2.6cm}
{\setlength{\unitlength}{0.3cm}
\begin{picture}(8,8)
\put(0,0){\includegraphics[width=2.4cm]{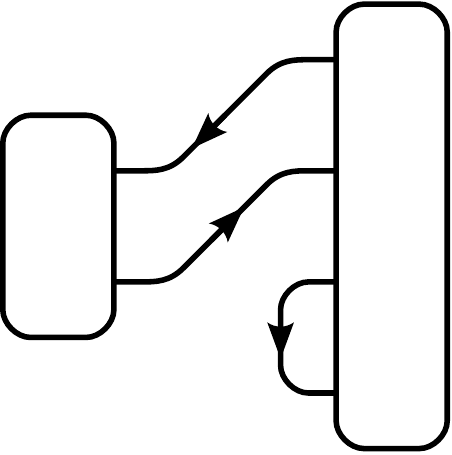}}
\put(0.5,3.6){\small $L_1$}
\put(6.5,3.6){\small $L_2$}
\end{picture}}
=
\parbox[c]{3.2cm}
{\setlength{\unitlength}{0.3cm}
\begin{picture}(10,8)
\put(0,0){\includegraphics[width=3cm]{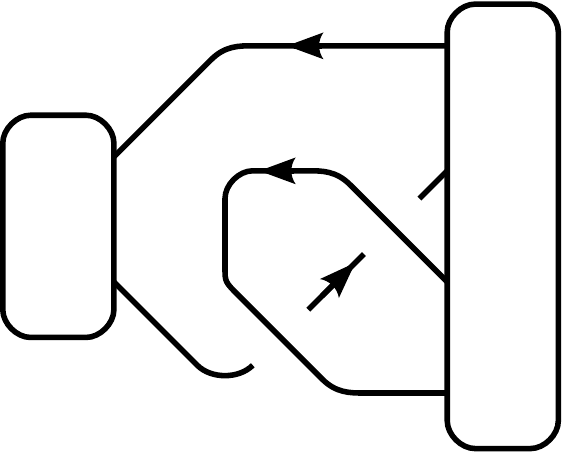}}
\put(0.5,3.6){\small $L_1$}
\put(8.5,3.6){\small $L_2$}
\end{picture}}\\
= x^{-1}t^{-1}
\parbox[c]{3.2cm}
{\setlength{\unitlength}{0.3cm}
\begin{picture}(10,8)
\put(0,0){\includegraphics[width=3cm]{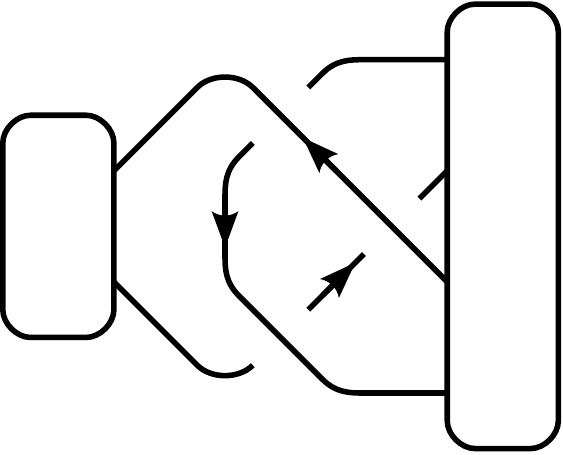}}
\put(0.5,3.6){\small $L_1$}
\put(8.5,3.6){\small $L_2$}
\end{picture}}
- x^{-1}t\,
\parbox[c]{3.2cm}
{\setlength{\unitlength}{0.3cm}
\begin{picture}(10,8)
\put(0,0){\includegraphics[width=3cm]{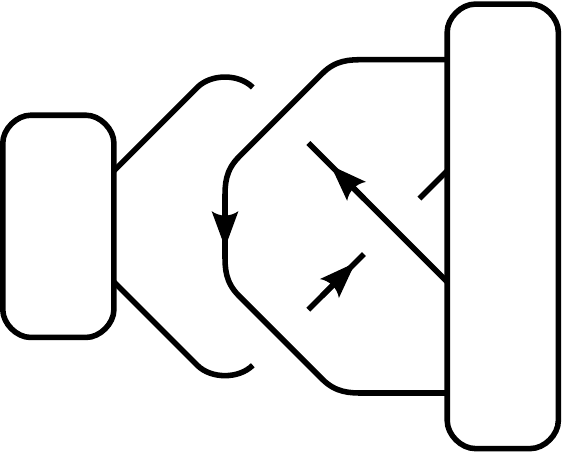}}
\put(0.5,3.6){\small $L_1$}
\put(8.5,3.6){\small $L_2$}
\end{picture}}\\
= x^{-1}t^{-1}
\parbox[c]{3.2cm}
{\setlength{\unitlength}{0.3cm}
\begin{picture}(10,8)
\put(0,0){\includegraphics[width=3cm]{ParWagV2F3c.pdf}}
\put(0.5,3.6){\small $L_1$}
\put(8.5,3.6){\small $L_2$}
\end{picture}}
- x^{-1}t\,
\parbox[c]{3.2cm}
{\setlength{\unitlength}{0.3cm}
\begin{picture}(10,8)
\put(0,0){\includegraphics[width=3cm]{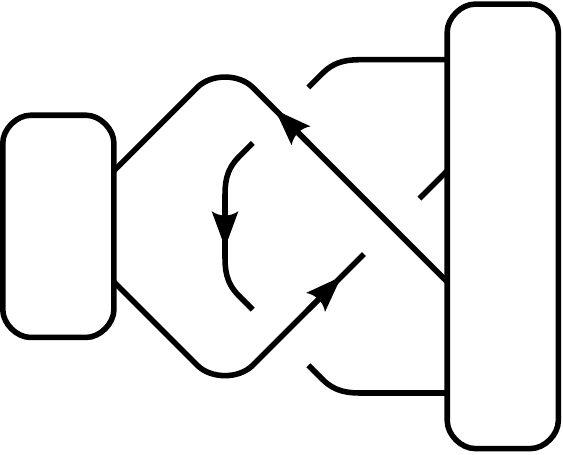}}
\put(0.5,3.6){\small $L_1$}
\put(8.5,3.6){\small $L_2$}
\end{picture}}\\
=
\parbox[c]{3.2cm}
{\setlength{\unitlength}{0.3cm}
\begin{picture}(10,8)
\put(0,0){\includegraphics[width=3cm]{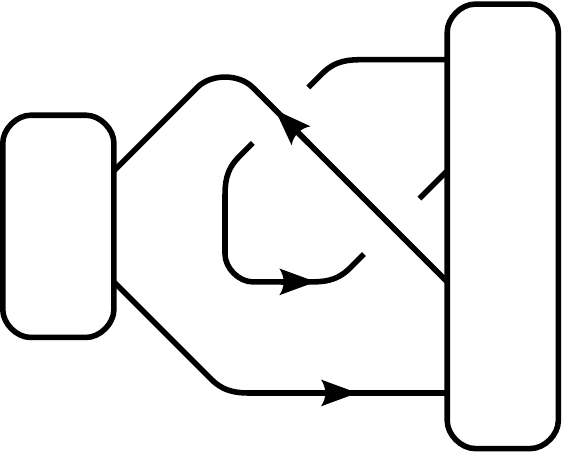}}
\put(0.5,3.6){\small $L_1$}
\put(8.5,3.6){\small $L_2$}
\end{picture}}
=
\parbox[c]{2.6cm}
{\setlength{\unitlength}{0.3cm}
\begin{picture}(8,8)
\put(0,0){\includegraphics[width=2.4cm]{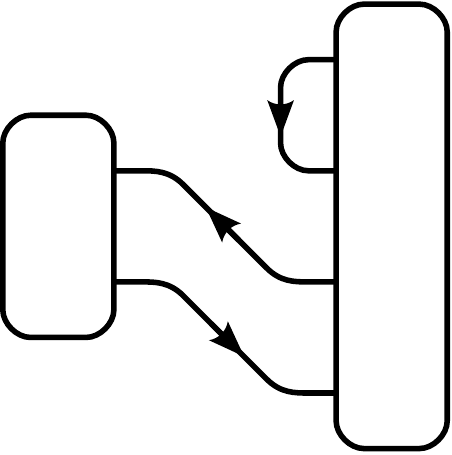}}
\put(0.5,3.6){\small $L_1$}
\put(6.5,3.6){\small $L_2$}
\end{picture}}
\end{gather*}
\centerline{{\bf Figure 2.4.} Graphical computation.}
\end{figure}

\bigskip\noindent
By Lemma 2.1, there is a well-defined operation $\star$ on $\SS(R,t,x)$ such that, if $L_1$ and $L_2$ are (classes of) singular links, then $L_1 \star L_2$ is represented by any connected sum of $L_1$ and $L_2$. Observe that $\SS(R,t,x)$ endowed with this operation is a graded algebra. Observe also that, thanks to (\ref{deuxoperations}), the operations $\cdot$ and $\star$ are linked by the following equality.

\bigskip\noindent
{\bf Lemma 2.2.}
{\it For all $u,v \in \SS(R,t,x)$ we have
\[
(t^{-1}-t) (u \star v) = x\,(u \cdot v)\,.
\]}
\qed

\bigskip\noindent
By Lemma 2.2, if $(t^{-1}-t)$ is invertible, then the operations $\star$ and $\cdot$ are essentially the same. Conversely, if $t = \pm 1$, then the operation $\cdot$ is trivial in the sense that the product of two elements is always zero, while $\SS(R,t,x)$ endowed with the operation $\star$ is a polynomial algebra. Indeed:

\bigskip\noindent
{\bf Theorem 2.3.}
{\it Let $R$ be a commutative ring, and let $t,x$ be two invertible elements in $R$. If $(t^{-1} -t -x)$ and $(t^{-1}-t+x)$ are invertible, then $\SS(R,t,x)$ endowed with the operation $\star$ is a polynomial algebra $R[X,Y]$ on two variables $X$ and $Y$, where $X$ and $Y$ are the (classes of) singular links represented in Figure 2.5.}

\begin{figure}[tbh]
\bigskip\bigskip\centerline{
\setlength{\unitlength}{0.4cm}
\begin{picture}(8,5)
\put(0,1){\includegraphics[width=3.2cm]{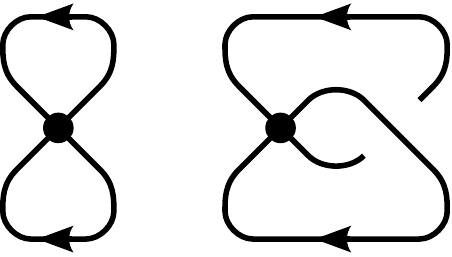}}
\put(0.6,0){\small $X$}
\put(5.6,0){\small $Y$}
\end{picture}} 

\smallskip
\centerline{{\bf Figure 2.5.} The links $X$ and $Y$.}
\end{figure}

\bigskip\noindent
{\bf Corollary 2.4.}
{\it Let $\K$ be a field.
\begin{itemize}
\item[(1)]
$\SS(\K(t,x),t,x)$ endowed with the operation $\cdot$ is a polynomial algebra $\K(t,x)[X,Y]$ on the two variables $X$ and $Y$.
\item[(2)]
$\SS(\K[x^{\pm 1}],1,x)$ endowed with the operation $\star$ is a polynomial algebra $\K[x^{\pm 1}][X,Y]$ on the two variables $X$ et $Y$.
\qed
\end{itemize}}

\bigskip\noindent
In order to prove Theorem 2.3, especially when the characteristic of $R$ is different from zero, we need to pass through  the ``ordered singular links'' and the ``ordered HOMFLY-PT skein module'' defined below. Note however that these two notions as well as Theorem 2.6 (see below) are interesting by themselves. 

\bigskip\noindent
Let $S(L)$ denote the set of singular points of a singular link $L$. An {\it ordered singular link} with $d$ singular points is a singular link $L$ with $d$ singular points together with a linear ordering on $S(L)$, that is, a one-to-one correspondence $o_L: S(L) \to \{1, \dots, d\}$. We denote by $\hat\LL_d$ the set of ordered singular links with $d$ singular points, and we set $\hat\LL = \sqcup_{d=0}^\infty \hat\LL_d$. Note that the symmetric group $\SSS_d$ acts naturally on $\hat\LL_d$, and $\hat\LL_d/\SSS_d = \LL_d$.

\bigskip\noindent
We denote by $R[\hat\LL_d]$ the free $R$-module freely generated by $\hat\LL_d$ for all $d \in \N$, and we denote by $R[\hat\LL]$ the free $R$-module freely generated by $\hat\LL$. We have $R[\hat\LL] = \oplus_{d=0}^\infty R[\hat\LL_d]$. For $d \in \N$, the action of $\SSS_d$ on $\hat\LL_d$ extends to an $R$-linear action on $R[\hat\LL_d]$, and we have $R[\hat\LL_d]/\SSS_d = R[\LL_d]$. 

\bigskip\noindent
Let $t,x$ be two invertible elements in $R$. The \emph{ordered HOMFLY-PT skein module} of ordered singular links (relative to the triple $(R,t,x)$), denoted by $\HSS(R,t,x)$, is defined to be the quotient of $R[\hat\LL]$ by the relations
\begin{equation}\label{skeinord}
x\, L_0 = t^{-1}\, L_+ - t\,L_-\,,
\end{equation}
for all triples of ordered singular links $(L_0,L_+,L_-)$ such that $L_0,L_+,L_-$ have the same singular points, ordered in the same way, and $L_0,L_+,L_-$ have the same diagram except in a disk where they are like in Figure 2.1. On the other hand, for $d \in \N$, we denote by $\HSS_d(R,t,x)$ the quotient of $R[\hat\LL_d]$ by the relations (\ref{skeinord}). Observe that
\[
\HSS(R,t,x) = \bigoplus_{d=0}^\infty \HSS_d(R,t,x)\,.
\]
The action of $\SSS_d$ on $R[\hat\LL_d]$ passes to the quotient, and thus induces an action on $\HSS_d(R,t,x)$. Moreover, $\HSS_d(R,t,x)/\SSS_d = \SS_d(R,t,x)$. For all $d \in \N$, we denote by $\pi_d : \HSS_d(R,t,x) \to \SS_d(R,t,x)$ the quotient map, and we denote by $\pi: \HSS(R,t,x) \to \SS(R,t,x)$ the $R$-linear epimorphism induced by the $\pi_d$'s.

\bigskip\noindent
Let $L_1$ and $L_2$ be two ordered singular links with $d_1$ and $d_2$ singular points, respectively. Let $L$ be a connected sum of $L_1$ and $L_2$, viewed as unordered singular links. By construction, we can assume that $S(L) = S(L_1) \sqcup S(L_2)$. Define a linear ordering on $S(L)$ by setting $o_L(p) = o_{L_1}(p)$ if $p \in S(L_1)$, and $o_L(p) = d_1+o_{L_2}(p)$ if $p \in S(L_2)$. The ordered singular link obtained in this way is called {\it connected sum} of the ordered singular links $L_1$ and $L_2$. The following is proved in the same way as Lemma 2.1. 

\bigskip\noindent
{\bf Lemma 2.5.}
{\it Two connected sums of two ordered singular links $L_1$ and $L_2$ represent the same element in $\HSS(R,t,x)$.}
\qed

\bigskip\noindent
By Lemma 2.5, there is a well-defined operation $\star$ on $\HSS(R,t,x)$ such that, if $L_1$ and $L_2$ are (classes of) ordered singular links, then $L_1 \star L_2$ is represented by any connected sum of $L_1$ and $L_2$. Observe that $\HSS(R,t,x)$ endowed with this operation is a graded algebra, and $\pi: \HSS(R,t,x) \to \SS(R,t,x)$ is an algebra epimorphism. Theorem 2.3 will be a consequence of the following.

\bigskip\noindent
{\bf Theorem 2.6.}
{\it Let $R$ be a commutative ring, and let $t,x$ be two invertible elements in $R$. If $(t^{-1} -t -x)$ and $(t^{-1}-t+x)$ are invertible in $R$, then $\HSS(R,t,x)$ endowed with the operation $\star$ is a free algebra $R \langle X,Y \rangle$ on $\{ X, Y\}$, where $X$ and $Y$ are the singular links represented in Figure~2.5.}


\section{The proofs}

From now on and throughout the whole section $R$ denotes a commutative ring and $t,x$ denote two invertible elements in $R$. 
The following lemma is known to experts, but we include a proof in order to reassure the reader.
Also, we insist on the fact that the proof and the conclusion of the lemma are valid for any $R$, even if its chracteristic is different from $0$.

\bigskip\noindent
{\bf Lemma 3.1.}
{\it The $R$-module $\SS_0(R,t,x) = \HSS_0(R,t,x)$ is isomorphic to $R$ and generated by $\bO$, where $\bO$ denotes the trivial knot.}

\bigskip\noindent
{\bf Proof.}
Let $\tilde x$ and $\tilde t$ be two variables. The classical study of the HOMFLY-PT polynomial on non-singular links (see \cite{HOMFLY, Jones2, Jones1, PrzTra1}) shows that there is a unique map $P: \LL_0 \to \Z[\tilde t^{\pm 1}, \tilde x^{\pm 1}]$ such that $P(\bO)=1$, and 
\[
\tilde x\, P(L_0) =\tilde t^{-1}\, P(L_+) -\tilde t\,P(L_-)\,,
\]
for all triples of links $(L_0,L_+,L_-)$ such that $L_0,L_+, L_-$ have the same diagram except in a disk where they are like in Figure 2.1. Let $\varphi : \Z[ \tilde t^{\pm 1}, \tilde x^{\pm 1}] \to R$ be the ring homomorphism which sends $\tilde t$ to $t$ and $\tilde x$ to $x$ (and $1$ to $1_R$). By the above, there is a well-defined $R$-homomorphism $f : \SS_0(R,t,x) \to R$ which sends $L$ to $\varphi(P(L))$ for all $L \in \LL_0$. It is easily seen that the $R$-homomorphism $g: R \to \SS_0(R,t,x)$ which sends $a \in R$ to $a\,\bO$ is a section of $f$, and that $\SS_0(R,t,x)$ is spanned by $\bO$, hence $f$ is an isomorphism.
\qed

\bigskip\noindent
With $\varepsilon = (\varepsilon_1, \dots, \varepsilon_d) \in \{0,1\}^d$ we associate the element $Z_\varepsilon$ in $\HSS_d(R,t,x)$ defined by
\[
Z_\varepsilon = Z_{\varepsilon_1} \star Z_{\varepsilon_2} \star \cdots \star Z_{\varepsilon_d}\,,
\]
where $Z_{\varepsilon_i} = X$ if $\varepsilon_i = 0$, and $Z_{\varepsilon_i} = Y$ if $\varepsilon_i =1$. 
For $d \ge 1$ we set $\BB_d=\{Z_\varepsilon \mid \varepsilon \in \{0,1\}^d \}$, and for $d=0$ we set $\BB_0=\{ \bO\}$. 

\bigskip\noindent
{\bf Lemma 3.2.}
{\it Let $d \in \N$. If $(t^{-1} -t -x)$ and $(t^{-1}-t+x)$ are invertible in $R$, then $\HSS_d(R,t,x)$ is spanned as an $R$-module by $\BB_d$.}

\bigskip\noindent
{\bf Proof.} 
We argue by induction on $d$. The case $d=0$ follows from Lemma 3.1, thus we can assume $d \ge 1$ plus the induction hypothesis. 

\bigskip\noindent
In our calculations we will need to consider the singular link $Y'$ represented in Figure 3.1. Observe that we have $t^{-1} Y - t Y'= x X$ in $\HSS_1(R,t,x)= \SS_1(R,t,x)$, hence $Y'$ lies in the $R$-submodule of $\HSS_1(R,t,x)$ spanned by $X$ and $Y$.

\begin{figure}[tbh]
\bigskip\bigskip\centerline{
\setlength{\unitlength}{0.4cm}
\begin{picture}(4,5)
\put(0,1){\includegraphics[width=1.6cm]{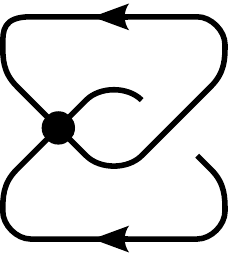}}
\put(1.5,0){\small $Y'$}
\end{picture}} 

\smallskip
\centerline{{\bf Figure 3.1.} The singular link $Y'$.}
\end{figure}

\bigskip\noindent
Let $U_d$ denote the submodule of $\HSS_d(R,t,x)$ spanned by $(X \star \HSS_{d-1} (R,t,x)) \cup (Y \star \HSS_{d-1} (R,t,x))$. By the induction hypothesis, in order to prove Lemma 3.2, it suffices to show that $\HSS_d(R,t,x) = U_d$. Let $L$ be an ordered singular link with $d$ singular points. We isolate the first singular point of $L$ like in Figure 3.2.

\begin{figure}[tbh]
\bigskip\bigskip\centerline{
\setlength{\unitlength}{0.3cm}
\begin{picture}(5,8)
\put(0,0){\includegraphics[width=1.5cm]{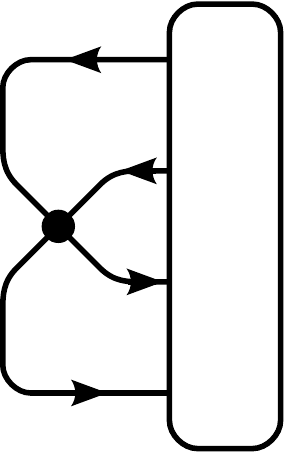}}
\put(0,3.5){\small $1$}
\put(3.4,3.5){\small $L_0$}
\end{picture}} 

\smallskip
\centerline{{\bf Figure 3.2.} First singular point of $L$.}
\end{figure}

\bigskip\noindent
We perform the graphical computations in $\HSS_d(R,t,x)$ shown in Figures 3.3 and 3.4.

\begin{figure}[tbh]
\begin{gather*}
\parbox[c]{1.7cm}
{\setlength{\unitlength}{0.3cm}
\begin{picture}(5,8)
\put(0,0){\includegraphics[width=1.5cm]{ParWagV2F6.pdf}}
\put(3.4,3.5){\small $L_0$}
\end{picture}}
=
\parbox[c]{3.2cm}
{\setlength{\unitlength}{0.3cm}
\begin{picture}(10,8)
\put(0,0){\includegraphics[width=3cm]{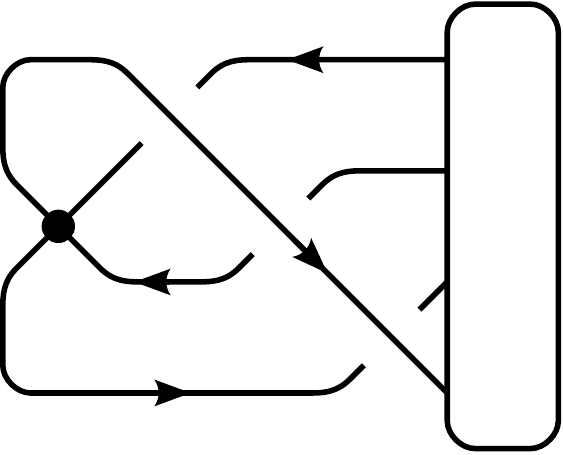}}
\put(8.4,3.5){\small $L_0$}
\end{picture}}\\
=t^{-2}
\parbox[c]{3.2cm}
{\setlength{\unitlength}{0.3cm}
\begin{picture}(10,8)
\put(0,0){\includegraphics[width=3cm]{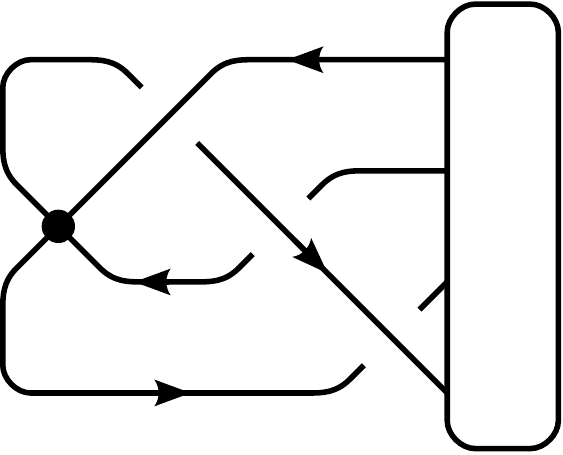}}
\put(8.4,3.5){\small $L_0$}
\end{picture}}
-xt^{-1}
\parbox[c]{3.2cm}
{\setlength{\unitlength}{0.3cm}
\begin{picture}(10,8)
\put(0,0){\includegraphics[width=3cm]{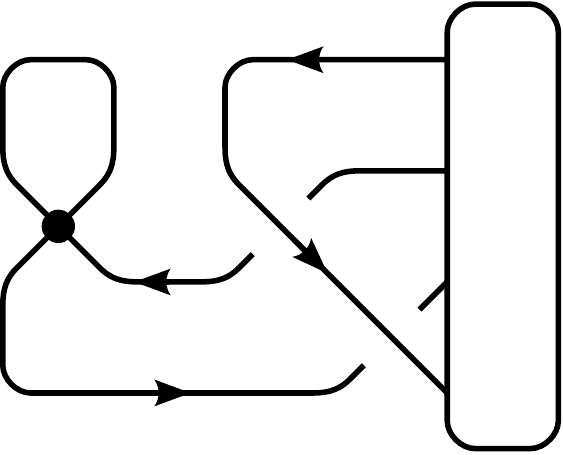}}
\put(8.4,3.5){\small $L_0$}
\end{picture}}\\
= t^{-4}
\parbox[c]{3.2cm}
{\setlength{\unitlength}{0.3cm}
\begin{picture}(10,8)
\put(0,0){\includegraphics[width=3cm]{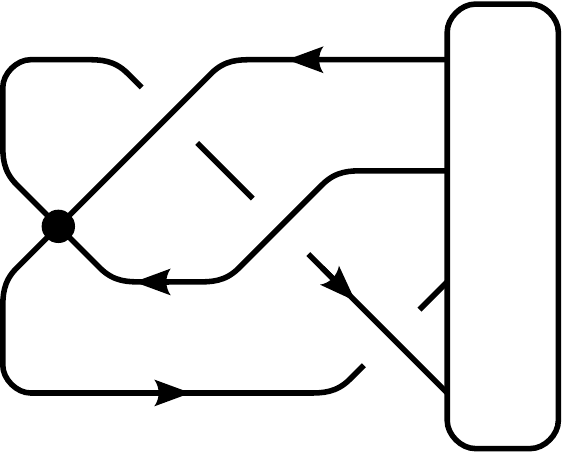}}
\put(8.4,3.5){\small $L_0$}
\end{picture}}
-xt^{-3}
\parbox[c]{3.2cm}
{\setlength{\unitlength}{0.3cm}
\begin{picture}(10,8)
\put(0,0){\includegraphics[width=3cm]{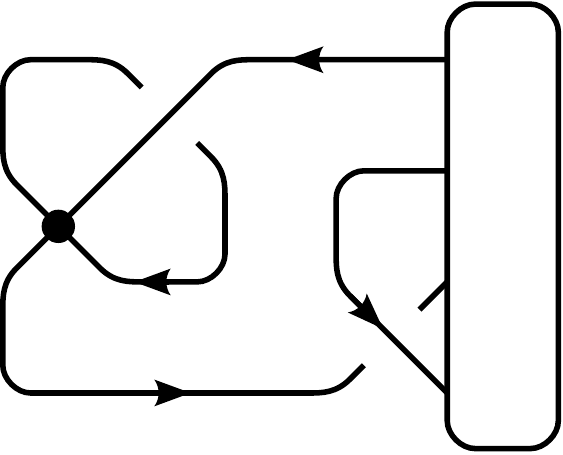}}
\put(8.4,3.5){\small $L_0$}
\end{picture}}
-xt^{-1}\, X \star
\parbox[c]{2cm}
{\setlength{\unitlength}{0.3cm}
\begin{picture}(6,8)
\put(0,0){\includegraphics[width=1.8cm]{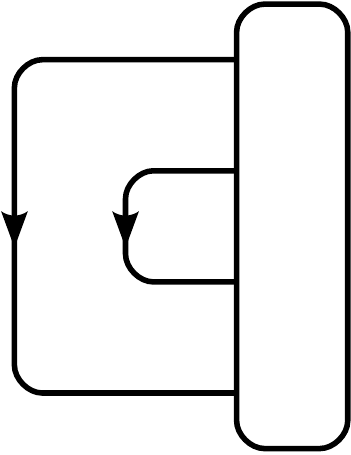}}
\put(4.4,3.5){\small $L_0$}
\end{picture}}\\
\equiv t^{-2}
\parbox[c]{3.2cm}
{\setlength{\unitlength}{0.3cm}
\begin{picture}(10,8)
\put(0,0){\includegraphics[width=3cm]{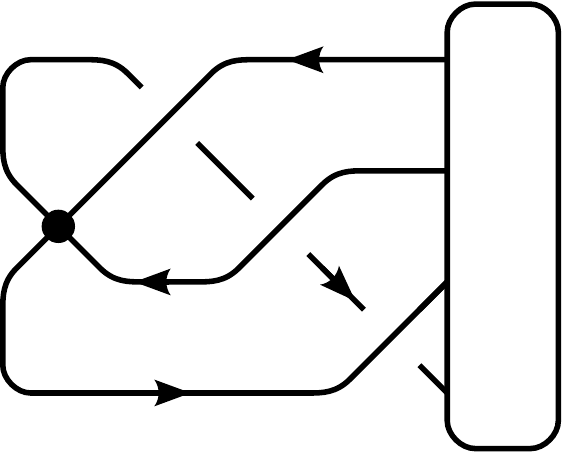}}
\put(8.4,3.5){\small $L_0$}
\end{picture}}
+ xt^{-3}
\parbox[c]{3.2cm}
{\setlength{\unitlength}{0.3cm}
\begin{picture}(10,8)
\put(0,0){\includegraphics[width=3cm]{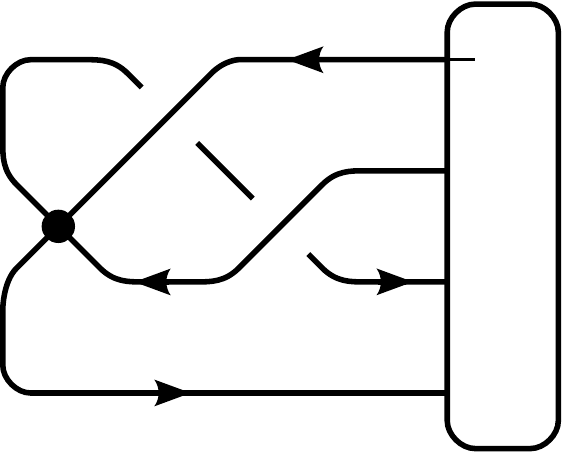}}
\put(8.4,3.5){\small $L_0$}
\end{picture}}
-xt^{-3}\, Y \star
\parbox[c]{2cm}
{\setlength{\unitlength}{0.3cm}
\begin{picture}(6,8)
\put(0,0){\includegraphics[width=1.8cm]{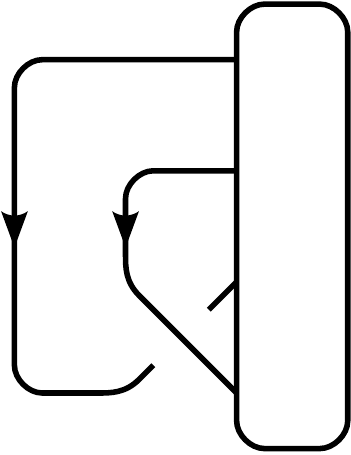}}
\put(4.4,3.5){\small $L_0$}
\end{picture}}
\ (\text{mod}\, U_d)\\
\equiv t^{-2}
\parbox[c]{1.7cm}
{\setlength{\unitlength}{0.3cm}
\begin{picture}(5,8)
\put(0,0){\includegraphics[width=1.5cm]{ParWagV2F6.pdf}}
\put(3.4,3.5){\small $L_0$}
\end{picture}}
+ x t^{-3}
\parbox[c]{2cm}
{\setlength{\unitlength}{0.3cm}
\begin{picture}(6,8)
\put(0,0){\includegraphics[width=1.8cm]{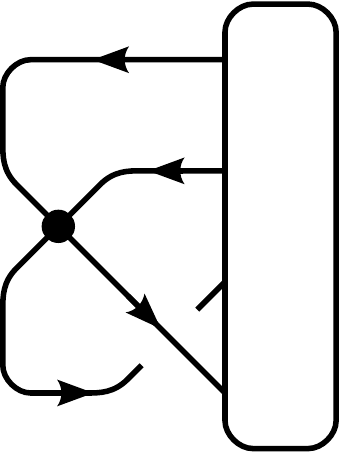}}
\put(4.4,3.5){\small $L_0$}
\end{picture}}
\ (\text{mod}\, U_d)
\end{gather*}
\centerline{{\bf Figure 3.3.} Graphical computation.}
\end{figure}

\begin{figure}[tbh]
\begin{gather*}
\parbox[c]{1.7cm}
{\setlength{\unitlength}{0.3cm}
\begin{picture}(5,8)
\put(0,0){\includegraphics[width=1.5cm]{ParWagV2F6.pdf}}
\put(3.4,3.5){\small $L_0$}
\end{picture}}
=
\parbox[c]{3.2cm}
{\setlength{\unitlength}{0.3cm}
\begin{picture}(10,8)
\put(0,0){\includegraphics[width=3cm]{ParWagV2F7g.pdf}}
\put(8.4,3.5){\small $L_0$}
\end{picture}}\\
= t^2
\parbox[c]{3.2cm}
{\setlength{\unitlength}{0.3cm}
\begin{picture}(10,8)
\put(0,0){\includegraphics[width=3cm]{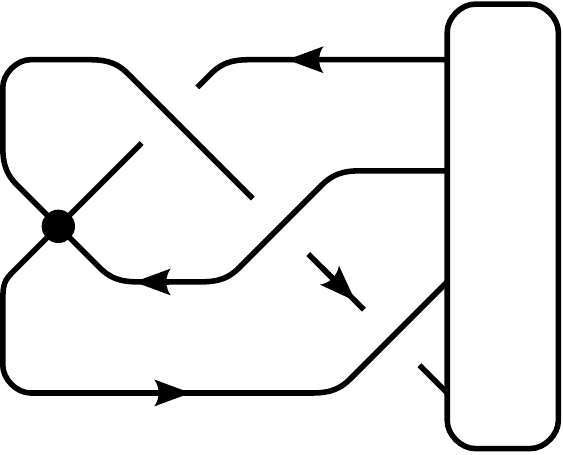}}
\put(8.4,3.5){\small $L_0$}
\end{picture}}
+xt
\parbox[c]{3.2cm}
{\setlength{\unitlength}{0.3cm}
\begin{picture}(10,8)
\put(0,0){\includegraphics[width=3cm]{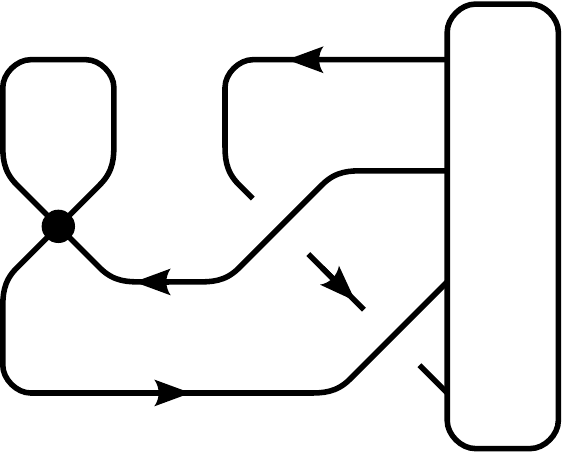}}
\put(8.4,3.5){\small $L_0$}
\end{picture}}\\
=t^4
\parbox[c]{3.2cm}
{\setlength{\unitlength}{0.3cm}
\begin{picture}(10,8)
\put(0,0){\includegraphics[width=3cm]{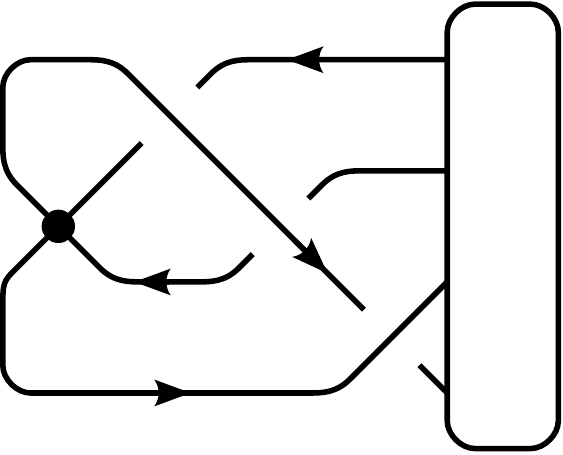}}
\put(8.4,3.5){\small $L_0$}
\end{picture}}
+xt^3
\parbox[c]{3.2cm}
{\setlength{\unitlength}{0.3cm}
\begin{picture}(10,8)
\put(0,0){\includegraphics[width=3cm]{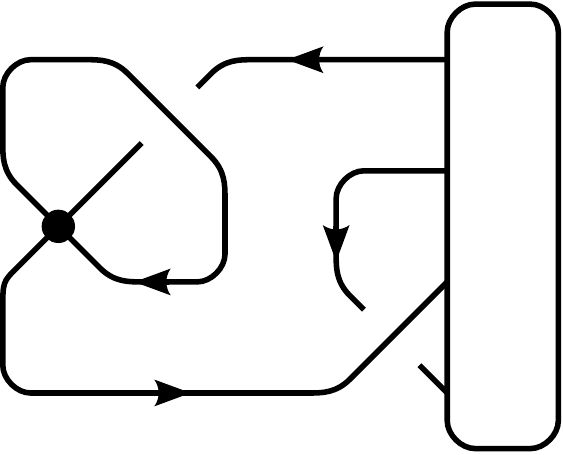}}
\put(8.4,3.5){\small $L_0$}
\end{picture}}
+xt\,X \star
\parbox[c]{2cm}
{\setlength{\unitlength}{0.3cm}
\begin{picture}(6,8)
\put(0,0){\includegraphics[width=1.8cm]{ParWagV2F7f.pdf}}
\put(4.4,3.5){\small $L_0$}
\end{picture}}\\
\equiv t^2
\parbox[c]{3.2cm}
{\setlength{\unitlength}{0.3cm}
\begin{picture}(10,8)
\put(0,0){\includegraphics[width=3cm]{ParWagV2F7a.pdf}}
\put(8.4,3.5){\small $L_0$}
\end{picture}}
-xt^3
\parbox[c]{3.2cm}
{\setlength{\unitlength}{0.3cm}
\begin{picture}(10,8)
\put(0,0){\includegraphics[width=3cm]{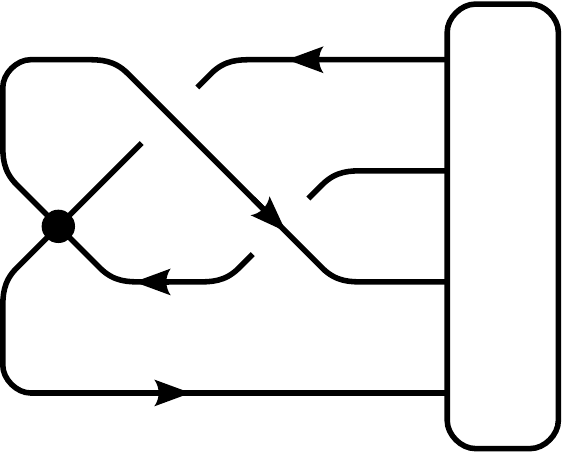}}
\put(8.4,3.5){\small $L_0$}
\end{picture}}
+xt^3\, Y'\star 
\parbox[c]{2cm}
{\setlength{\unitlength}{0.3cm}
\begin{picture}(6,8)
\put(0,0){\includegraphics[width=1.8cm]{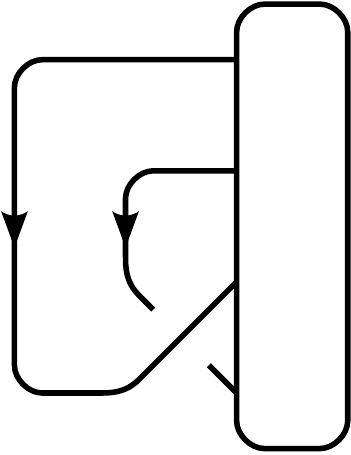}}
\put(4.4,3.5){\small $L_0$}
\end{picture}}
\ (\text{mod}\, U_d)\\
\equiv t^2
\parbox[c]{1.7cm}
{\setlength{\unitlength}{0.3cm}
\begin{picture}(5,8)
\put(0,0){\includegraphics[width=1.5cm]{ParWagV2F6.pdf}}
\put(3.4,3.5){\small $L_0$}
\end{picture}}
-xt^3
\parbox[c]{2cm}
{\setlength{\unitlength}{0.3cm}
\begin{picture}(6,8)
\put(0,0){\includegraphics[width=1.8cm]{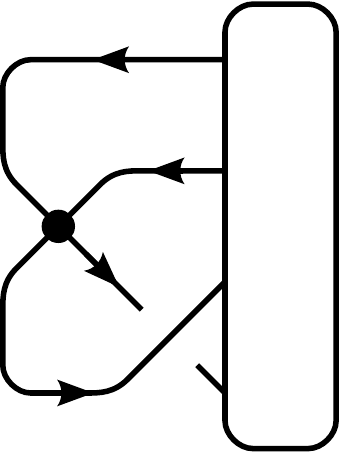}}
\put(4.4,3.5){\small $L_0$}
\end{picture}}
\ (\text{mod}\, U_d)
\end{gather*}
\centerline{{\bf Figure 3.4.} Graphical computation.}
\end{figure}

\bigskip\noindent
Summing $t^2$ times the congruence in Figure 3.3 with $t^{-2}$ times the congruence in Figure 3.4 we obtain the congruence in Figure 3.5. 

\begin{figure}[tbh]
\begin{align*}
(t^2+t^{-2} -2) 
\parbox[c]{1.7cm}
{\setlength{\unitlength}{0.3cm}
\begin{picture}(5,8)
\put(0,0){\includegraphics[width=1.5cm]{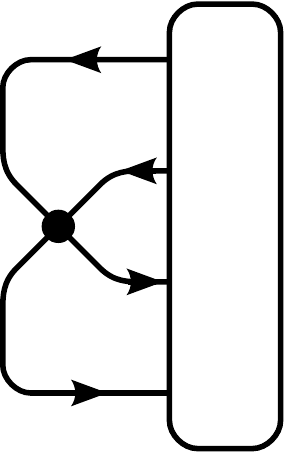}}
\put(3.4,3.5){\small $L_0$}
\end{picture}}
& \equiv xt^{-1} 
\parbox[c]{2cm}
{\setlength{\unitlength}{0.3cm}
\begin{picture}(6,8)
\put(0,0){\includegraphics[width=1.8cm]{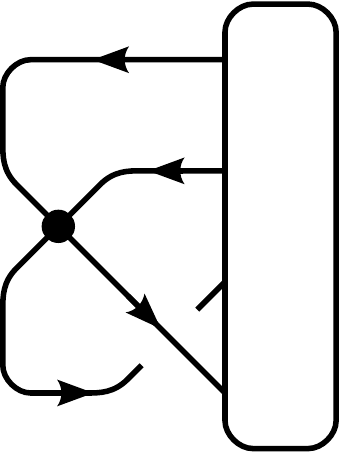}}
\put(4.4,3.5){\small $L_0$}
\end{picture}}
-xt
\parbox[c]{2cm}
{\setlength{\unitlength}{0.3cm}
\begin{picture}(6,8)
\put(0,0){\includegraphics[width=1.8cm]{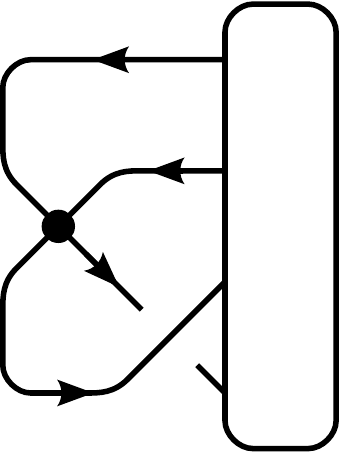}}
\put(4.4,3.5){\small $L_0$}
\end{picture}}
\ (\text{mod}\, U_d)\\
& \equiv x^2
\parbox[c]{1.7cm}
{\setlength{\unitlength}{0.3cm}
\begin{picture}(5,8)
\put(0,0){\includegraphics[width=1.5cm]{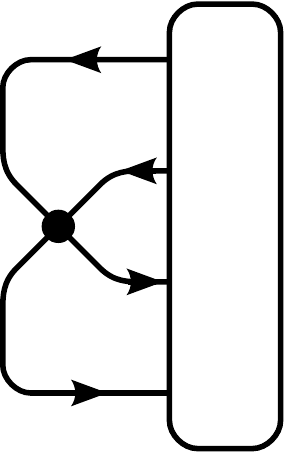}}
\put(3.4,3.5){\small $L_0$}
\end{picture}}
\ (\text{mod}\, U_d)
\end{align*}
\centerline{{\bf Figure 3.5.} Graphical computation.}
\end{figure}

\bigskip\noindent
We conclude that $(t^2+t^{-2}-2 -x^2)\,L \in U_d$, thus $L \in U_d$ since $(t^2+t^{-2}-2-x^2) = (t^{-1}-t -x)(t^{-1} -t+x)$ is invertible.
\qed

\bigskip\noindent
{\bf Lemma 3.3.}
{\it Let $d \in \N$. If $(t^{-1} -t -x)$ and $(t^{-1}-t+x)$ are invertible in $R$, then $\BB_d$ is $R$-linearly free in $\HSS_d(R,t,x)$.}

\bigskip\noindent
{\bf Proof.}
We argue by induction on $d$. By Lemma 3.1, $\{\bO\}$ is linearly free in $\HSS_0(R,t,x)$. So, we can assume $d \ge 1$ plus the induction hypothesis. Let $L\ \in \hat \LL_d$. Let $p$ be the first singular point of $L$. We denote by $g_{0,d}(L)$ the ordered singular link in $\hat\LL_{d-1}$ obtained from $L$ by replacing a neighborhood of $p$ by a non-crossing, like in Figure 3.6 (left hand side), and we denote by $g_{1,d}(L)$ the ordered singular link in $\hat\LL_{d-1}$ obtained from $L$ by replacing a neighborhood of $p$ by a negative crossing, like in Figure 3.6 (right hand side). For $L' \in \{g_{0,d}(L), g_{1,d}(L)\}$, the singular points of $L'$ are ordered by setting $o_{L'}(q) = o_L(q)-1$ for all $q \in S(L')$. It is easily checked using the singular Reidemeister moves (see \cite{Kauff1}) that the map $g_{i,d} : \hat\LL_d \to \hat\LL_{d-1}$ is well-defined. We extend $g_{i,d}$ linearly, $g_{i,d} : R[\hat\LL_d] \to R[\hat\LL_{d-1}]$, and we observe that $g_{i,d}$ induces a $R$-linear map $f_{i,d}: \HSS_d(R,t,x) \to \HSS_{d-1}(R,t,x)$. Note that, if $u \in \HSS_1(R,t,x)$ and $v \in \HSS_{d-1}(R,t,x)$, then 
\[
f_{i,d}(u \star v) = f_{i,1}(u) \star v\,.
\]
Moreover,
\[
f_{0,1}(X) = x^{-1}(t^{-1}-t)\, \bO\,,\ f_{0,1}(Y) = \bO\,,\ f_{1,1}(X) = \bO\,, \text{ and } f_{1,1}(Y)= x^{-1}(t^{-1}-t)\, \bO\,.
\]

\begin{figure}[tbh]
\bigskip\bigskip\centerline{
\setlength{\unitlength}{0.4cm}
\begin{picture}(18,2)
\put(0,0){\includegraphics[width=7.2cm]{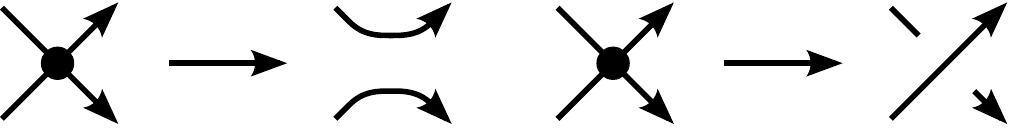}}
\put(0.6,0){\small $p$}
\put(3,1.6){\small $g_{0,d}$}
\put(10.6,0){\small $p$}
\put(13,1.6){\small $g_{1,d}$}
\end{picture}} 

\smallskip
\centerline{{\bf Figure 3.6.} The maps $g_{0,d}$ and $g_{1,d}$.}
\end{figure}

\bigskip\noindent
Assume we have an equality
\[
\sum_{\varepsilon \in \{0,1\}^d} a_\varepsilon Z_\varepsilon = 0
\]
in $\HSS_d(R,t,x)$, where $a_\varepsilon \in R$ for all $\varepsilon \in \{0,1\}^d$. Applying $f_{0,d}$ to this equality we get 
\[
0 = \sum_{\mu \in \{0,1\}^{d-1}} \big( a_{(0,\mu)} f_{0,1}(X) + a_{(1,\mu)} f_{0,1}(Y)) \star Z_\mu
= \sum_{\mu \in \{0,1\}^{d-1}} \big(x^{-1}(t^{-1}-t) a_{(0,\mu)} + a_{(1,\mu)}) Z_\mu\,.
\]
By the induction hypothesis, this last equality implies
\begin{equation}\label{Egalite1}
x^{-1}(t^{-1}-t) a_{(0,\mu)} + a_{(1,\mu)}=0
\end{equation}
for all $\mu \in \{0,1\}^{d-1}$. Similarly, applying $f_{1,d}$ we get 
\[
0 = \sum_{\mu \in \{0,1\}^{d-1}} \big( a_{(0,\mu)} f_{1,1}(X) + a_{(1,\mu)} f_{1,1}(Y)) \star Z_\mu
= \sum_{\mu \in \{0,1\}^{d-1}} \big(a_{(0,\mu)} +x^{-1}(t^{-1}-t)  a_{(1,\mu)}) Z_\mu\,,
\]
hence 
\begin{equation}\label{Egalite2}
a_{(0,\mu)} +x^{-1}(t^{-1}-t)  a_{(1,\mu)}=0
\end{equation}
for all $\mu \in \{0,1\}^{d-1}$. Equalities (\ref{Egalite1}) and (\ref{Egalite2}) form the following system of linear equations. 
\[ 
\left\{\begin{array}{rrcl}
(t^{-1}-t)\, a_{(0,\mu)} & +x\, a_{(1,\mu)}& = & 0 \\
x\, a_{(0,\mu)} & + (t^{-1}-t)\,  a_{(1,\mu)} & = & 0
\end{array} \right.
 \]
The determinant of this system is $(t^{-1}-t-x)(t^{-1}-t+x)$, which is invertible, thus $a_{(0,\mu)} = a_{(1,\mu)} = 0$ for all $\mu \in \{0,1\}^{d-1}$. This shows that $\BB_d$ is linearly free in $\HSS_d(R,t,x)$.
\qed

\bigskip\noindent
{\bf Proof of Theorem 2.6.}
The equality
\[
\HSS(R,t,x) = \bigoplus_{d=0}^\infty \HSS_d(R,t,x)
\]
equips $(\HSS(R,t,x), \star)$ with a structure of graded algebra. On the other hand, by Lemmas 3.2 and 3.3, for all $d \in \N$, $\HSS_d(R,t,x)$ is the free $R$-module freely generated by $\BB_d$. It follows that $\HSS(R,t,x)$ is the free algebra $R \langle X,Y \rangle$ on $\{ X, Y\}$.
\qed

\bigskip\noindent
{\bf Proof of Theorem 2.3.}
The symmetric group $\SSS_d$ acts on $\{0,1\}^d$ by permutations of the coordinates, and we have $w\, Z_\varepsilon = Z_{w \varepsilon}$ for all $w \in \SSS_d$ and $\varepsilon \in \{0,1\}^d$. The fact that the action of $\SSS_d$ permutes the elements of the $R$-basis $\BB_d$ of $\HSS_d(R,t,x)$ implies that $\HSS_d(R,t,x)/\SSS_d$ is the free $R$-module freely generated by the orbits of $\SSS_d$ in the basis. In other words, $\SS_d(R,t,x) = \HSS_d(R,t,x)/\SSS_d$ is the free $R$-module freely generated by $\{ X^{\star i} \star Y^{\star j} \mid i+j=d\}$. We conclude that $\SS(R,t,x) = \oplus_{d=0}^\infty \SS_d(R,t,x)$ is the polynomial algebra $R[X,Y]$ on the variables $X,Y$.
\qed



\bigskip\bigskip\noindent
{\bf Luis Paris,}

\smallskip\noindent 
Université de Bourgogne, Institut de Mathématiques de Bourgogne, UMR 5584 du CNRS, B.P. 47870, 21078 Dijon cedex, France.

\smallskip\noindent
E-mail: {\tt lparis@u-bourgogne.fr}

\bigskip\noindent
{\bf Emmanuel Wagner,}

\smallskip\noindent 
Université de Bourgogne, Institut de Mathématiques de Bourgogne, UMR 5584 du CNRS, B.P. 47870, 21078 Dijon cedex, France.

\smallskip\noindent
E-mail: {\tt emmanuel.wagner@u-bourgogne.fr}

\end{document}